\documentclass[11pt,reqno]{amsart}
\usepackage{amssymb,amsmath,mathrsfs}
\usepackage{enumerate,tikz}
\usepackage{hyperref}
\usepackage{forest}

\usetikzlibrary{decorations.pathreplacing}

\newtheorem{thm}{Theorem}[section]
\newtheorem{prop}[thm]{Proposition}

\theoremstyle{remark}
\newtheorem*{remark}{Remark}

\numberwithin{equation}{section}

\def\U{\mathsf{U}}
\def\D{\mathsf{D}}
\def\B{\mathcal{B}}
\def\C{\mathcal{C}}
\def\DP{\mathcal{D}}
\def\Pol{\mathcal{P}}
\def\T{\mathcal{T}}

\def\abs#1{\lvert#1\rvert}
\def\oeis#1{\cite[#1]{Sloane}}

\definecolor{red1}{HTML}{B02400}
\def\bcol#1{\textcolor{cyan}{#1}}
\def\rcol#1{\textcolor{red1}{#1}}
\def\reddot{\tikz{\draw[red1, fill] circle(0.01)}}

\DeclareMathOperator{\val}{val}

\newcommand\dyckpath[4]{
\def\diam{0.08}
\def\sf{#4}
  \begin{tikzpicture}[scale=\sf]
    \draw[help lines] (#1) -- ++(#2*2,0);
    \draw[line width=1pt] (#1) foreach \dir in {#3}{ -- ++(\dir*90-45:1.41)};
    \draw[fill] (#1) circle (\diam);
    \draw[fill] (#1) foreach \dir in {#3}{ ++(\dir*90-45:1.41) circle (\diam)};
   \end{tikzpicture}
}

\begin{document}
\title[Colored compositions, Dyck paths, and polygon partitions]{Bijections between colored compositions, Dyck paths, and polygon partitions}
\author[Juan Gil]{Juan B. Gil}
\address{Penn State Altoona\\ 3000 Ivyside Park\\ Altoona, PA 16601}
\email{jgil@psu.edu}

\author[Emma Hoover]{Emma G. Hoover}
\email{egh5128@psu.edu}

\author[Jessica Shearer]{Jessica A. Shearer}
\email{jas8642@psu.edu}

\begin{abstract}
In this paper, we give part-preserving bijections between three fundamental families of objects that serve as natural framework for many problems in enumerative combinatorics. Specifically, we consider compositions, Dyck paths, and partitions of a convex polygon, and identify suitable building blocks that are then appropriately decorated to achieve matching cardinalities. Our bijections are constructive and apply for the general case where the building blocks are allowed to come in different colors.
\end{abstract}

\maketitle

%%%%%%%%%%%%%%%%%%%%%%%%%%%%%%%%%%%%%%%%%%%%%%%%%%%%%
\section{Introduction}
Motivated by the approach suggested by Birmajer, Gil, and Weiner \cite{BGW19}, we use known enumeration formulas in terms of partial Bell polynomials to establish connections between colored compositions, colored Dyck paths, and colored partitions of convex polygons. 

Let us start by recalling some basic definitions and known results.

\subsection*{Compositions}
A {\em composition} of a positive integer $n$ is an ordered $k$-tuple $(j_1, \dotsc, j_k)$, for $k \ge 1$, of positive integers (parts) such that $j_1 + \dotsb + j_k = n$. For example, $(1, 3, 3, 2)$ is a composition of $n=9$ with $k=4$ parts. 

There are $\binom{n-1}{k-1}$ compositions of $n$ with $k$ parts. Moreover, given a sequence of nonnegative integers $\gamma=(\gamma_j)_{j\in\mathbb{N}}$, a $\gamma$-color composition of $n$ is a composition such that part $j$ can take on $\gamma_j$ colors. If $\gamma_j = 0$, then $j$ is not admissible as part of the composition. We let $c_{n,k}(\gamma)$ denote the number of $\gamma$-color compositions of $n$ with exactly $k$ parts.

\begin{prop}[Hoggatt \& Lind \cite{HoLi68}]
\[ c_{n,k}(\gamma) = \frac{k!}{n!} B_{n, k} (1!\gamma_1, 2! \gamma_2, \dots) \;\text{ for } n\ge 1, \]
where $B_{n,k}$ is the $(n,k)$-th partial Bell polynomial.\footnote{For the definition and properties of these polynomials, see the book by Comtet \cite[Sec.~3.3]{Comtet}.}
\end{prop}

\subsection*{Dyck paths}
A {\em Dyck path} of semilength $n$ (or $n$-Dyck path) is a lattice path from $(0,0)$ to $(2n,0)$, consisting of steps $(1,1)\leadsto\U$ and $(1,-1)\leadsto\D$, never going below the $x$-axis. For example, the lattice path
\begin{center}
\begin{tikzpicture}
\node at (0,0) {\dyckpath{0,0}{9}{1,1,1,1,0,0,0,1,0,0,1,1,0,1,1,0,0,0}{0.35}};
\end{tikzpicture}
\end{center}
is a Dyck path of semilength $n=9$ with 4 peaks.

Note that every Dyck path can be constructed using primitive blocks of the form 
\[ \U, \U\D, \U\D^2, \U\D^3, \U\D^4,\dots \]
or alternatively, using primitive blocks of the form 
\[ \D, \U\D, \U^2\D, \U^3\D, \U^4\D,\dots \]
For instance, our sample Dyck path could be written as $\U\U\U(\U\D^3)(\U\D^2)\U(\U\D)\U(\U\D^3)$ or as $(\U^4\D)\D\D(\U\D)\D(\U^2\D)(\U^2\D)\D\D$, depending on the chosen set of building blocks.

There are $\frac{1}{n} \binom{n}{k-1}\binom{n}{k}$ $n$-Dyck paths with exactly $k$ peaks. Given a sequence of nonnegative integers $\gamma=(\gamma_j)_{j\in\mathbb{N}}$, we let $d_{n,k}(\gamma)$ denote the number of $n$-Dyck paths with $k$ peaks such that the primitive component of length $j+1$ (e.g.\ $\U\D^j$) may be colored in $\gamma_j$ different ways.

\begin{prop}[Mansour \& Sun \cite{MaSu08}]
\[  d_{n,k}(\gamma) = \frac{1}{(n-k+1)!} B_{n,k}(1!\gamma_1, 2!\gamma_2, \dots) \text{ for } n\ge 1. \]
\end{prop}

\subsection*{Polygon partitions}
Finally, we consider the set of rooted convex polygons partitioned by noncrossing diagonals. For example,

\begin{center}
\begin{tikzpicture}
\node[regular polygon, regular polygon sides=11, minimum size=16ex] at (0,0) (A) {};
\draw[gray] (A.corner 1) -- (A.corner 3);
\draw[gray] (A.corner 3) -- (A.corner 6);
\draw[gray] (A.corner 6) -- (A.corner 10);
\node[regular polygon, regular polygon sides=11, minimum size=16.2ex, draw, thick] at (0,0) (B) {};
\draw[thick,gray,dash pattern=on 3pt off 1pt] (B.corner 6) -- (B.corner 7);
\end{tikzpicture}
\end{center}
is a partition of a rooted\footnote{Note that the bottom side is designated as the base/root.} undecagon into two pentagons, a quadrilateral, and a triangle.

Given a sequence of nonnegative integers $\gamma=(\gamma_j)_{j\in\mathbb{N}}$, we let $p_{n,k}(\gamma)$ be the number of colored partitions of a rooted convex $(n+2)$-gon made by $k-1$ noncrossing diagonals into $k$ polygons such that each $(j+2)$-gon may be colored in $\gamma_j$ different ways.
\begin{prop}[Birmajer, Gil \& Weiner \cite{BGW17a}]
\[ p_{n,k}(\gamma) =\frac1{n+1}\binom{n+k}{k}\frac{k!}{n!}B_{n,k}(1!\gamma_1,2!\gamma_2,\dots) \;\text{ for } n\ge 1. \]
\end{prop}

\bigskip
It is clear that the aforementioned sets are not equinumerous, but with appropriate markings, they can be bijectively connected. In fact, based on the above three propositions, it is straightforward to see that for any coloring sequence $\gamma=(\gamma_j)_{j\in\mathbb{N}}$, we have:
\begin{equation}\label{eq:mainIdentities}
\begin{aligned}
 \binom{n}{k-1} c_{n,k}(\gamma) &= k\,d_{n,k}(\gamma), \\[5pt]
 \binom{n+k}{k} d_{n,k}(\gamma) &= \binom{n+1}{k} p_{n,k}(\gamma), \\[5pt]
 (n+1) p_{n,k}(\gamma) &= \binom{n+k}{k} c_{n,k}(\gamma).
\end{aligned}
\end{equation}
The goal of this paper is to provide constructive combinatorial proofs for these identities.

%%%%%%%%%%%%%%%%%%%%%%%%%%%%%%%%%%%%%%%%%%%%%%%%%%%%%
\section{Bijection between compositions and Dyck paths}
\label{sec:Comp2Dyck}

In this section, we give a bijective proof of the following result.

\begin{thm}
Let $1\le k\le n$. For any coloring sequence $\gamma=(\gamma_j)_{j\in\mathbb{N}}$,
\begin{equation*}
\binom{n}{k-1}c_{n,k}(\gamma) = k\,d_{n,k}(\gamma).
\end{equation*}
\end{thm}

\medskip
We first assume $\gamma=(1,1,1,\dots)$ and represent every composition of $n$ as a tiling of an $n\times 1$ rectangular board with $n$ unit squares. For example,
\begin{center}
\begin{tikzpicture}[scale=0.5]
\draw [step=1,thin,gray!30] (0,0) grid (9,1);
\draw [thick] (0,0) rectangle (9,1); 
\foreach \x in {1,4,7}{\draw[thick] (\x,0) -- (\x,1);} 
\node [left=5pt] at (0,0.5) {(1, 3, 3, 2) \; $\longleftrightarrow$}; 
\end{tikzpicture}
\end{center}

Moreover, we can represent every lattice path with steps $\U$ and $\D$ as a word over the alphabet $\{\U,\D\}$ with valuations $\val(\U)=1$ and $\val(\D)=-1$. A Dyck word $w$ is then a word over that alphabet with $\val(w)=0$ and such that, if $w=uv$, then $\val(u)\ge 0$.

Let $\C_{n,k}^{(k-1)}$ be the set of compositions of $n$ with $k$ parts and such that $k-1$ of the $n$ unit squares are marked. Moreover, let $\DP_{n,k}^{(\textup{peak})}$ be the set of $n$-Dyck paths with $k$ peaks and such that one peak is marked. 

We now proceed to give a bijective map $\varphi_{dc}: \DP_{n,k}^{(\textup{peak})}\to \C_{n,k}^{(k-1)}$, illustrating each step of our construction with an example.
\begin{enumerate}[\quad$\circ$]
\item Given a Dyck path $P$ with one peak marked, write it as a word $w_P$ and extend it by adding an extra $\U$ to the beginning $\leadsto w=\U w_P$. Note that $\val(w)=1$, and if $w=uv$ with a nonempty word $u$, then $\val(u)>0$. 
\begin{center}
\begin{tikzpicture}
\node at (0,0) {\dyckpath{0,0}{9}{1,1,1,1,0,0,0,1,0,0,1,1,0,1,1,0,0,0}{0.4}};
\draw[fill,red1] (1.32,0) circle(0.06);
\node[below=1] at (0,-1) {$\leadsto w = \bcol{\U}\U\U\U\U\D\D\D\U\D\D\U\U\D\U\U\D\D\D$};
\end{tikzpicture}
\end{center}

\item Split the word $w$ at each valley (occurrence of $\D\U$) into subwords $w_1,\dots,w_k$ that start with a string of $\U$'s and end with a string of $\D$'s, and mark the subword that corresponds to the marked peak of the path.
\begin{center}
\tikz[baseline=0]{\node at (0,0){$(\U\U\U\U\U\D\D\D)(\U\D\D)(\U\U\D)(\U\U\D\D\D)$};
\draw[red1,thick] (0.5,0.25)--(1.4,0.25);}
\end{center}
\item If needed, rotate the subwords $w_1,\dots,w_k$, repeatedly moving the last subword to the front of the word until the marked subword is the left most subword. 
\begin{center}
\tikz[baseline=0]{\node at (0,0){$(\U\U\U\U\U\D\D\D)(\U\D\D)(\U\U\D)(\U\U\D\D\D)$};
\draw[red1,thick] (0.5,0.25)--(1.4,0.25);}\\
$\downarrow$ \\
\tikz[baseline=0]{\node at (0,0){$(\U\U\D\D\D)(\U\U\U\U\U\D\D\D)(\U\D\D)(\U\U\D)$};
\draw[red1,thick] (2.15,0.25)--(3.05,0.25);}\\
$\downarrow$ \\
\tikz[baseline=0]{\node at (0,0){$(\U\U\D)(\U\U\D\D\D)(\U\U\U\U\U\D\D\D)(\U\D\D)$};
\draw[red1,thick] (-3.08,0.25)--(-2.17,0.25);}
\end{center}
Note that each of the rotated words will have a left factor with valuation $0$. 
\item Let $\hat w$ denote the (rotated) word starting with the marked subword. The $k$ parts of the composition associated to $P$ are the number of $\D$'s in each of the $k$ subwords of $\hat w$.
\smallskip
\begin{center}
\begin{tikzpicture}[scale=0.6]
\draw [step=1,thin,gray!40] (0,0) grid (9,1);
\draw [thick] (0,0) rectangle (9,1); 
\foreach \x in {1,4,7}{\draw[thick] (\x,0) -- (\x,1);}
\end{tikzpicture}
\end{center}
\item Next, we identify the first $k-1$ peaks (occurrences of $\U\D$) in $\hat w$ and mark the positions of the corresponding $\U$-steps. If $k=1$, we skip this step.
\begin{center}
\tikz[baseline=0]{\node at (0,0){$(\U\U\D)(\U\U\D\D\D)(\U\U\U\U\U\D\D\D)(\U\D\D)$};
\foreach \x in {-2.72,-1.62,0.84}{\draw[red1,very thick] (\x,-0.18) -- (\x+0.2,-0.18);}}
\end{center}
\item Finally, identifying the first $n$ $\U$-steps of $\hat w$ with the $n$ unit squares of the associated composition, we place $k-1$ dots at the  unit squares that correspond to the positions of the marked peaks in $\hat w$.
\smallskip
\begin{center}
\begin{tikzpicture}[scale=0.6]
\draw [step=1,thin,gray!40] (0,0) grid (9,1);
\draw [thick] (0,0) rectangle (9,1); 
\foreach \x in {1,4,7}{\draw[thick] (\x,0) -- (\x,1);}
\foreach \m in {2,4,9}{\draw[fill,red1] (\m-0.5,0.5) circle(0.1);}
\end{tikzpicture}
\end{center}
The resulting marked composition is denoted by $\varphi_{dc}(P)$.
\end{enumerate}

Going back from a marked composition $C$ to a word $\hat w$ over the alphabet $\{\U,\D\}$ is fairly simple. The length of each tile in $C$ gives the distribution of the $\D$'s, and the location of the dots gives the placement of the $\U$'s (with an extra $\U$ added to the last group of $\D$'s). By construction, $\val(\hat w)=1$. 

Now, split $\hat w$ at its valleys to obtain subwords $\hat w_1,\dots,\hat w_k$ that start with a string of $\U$'s and end with a string of $\D$'s. Mark the subword $\hat w_1$. If $\hat w$ has a left factor with valuation $0$, we rotate the subwords repeatedly moving the first subword to the end of the word until the rotated subword has no left factor with valuation $0$. Once that point is reached, we remove the first $\U$ from the rotated word and obtain a marked Dyck word, $\varphi_{dc}^{-1}(C)$. Note that the marked peak is the one inside the subword $\hat w_1$.  

We finish by observing that, under our bijection, each part $j$ of a given composition corresponds to a maximal descent $\U\D^j$ of the associated Dyck path, so our proof also works for an arbitrary coloring sequence $\gamma = (\gamma_1, \gamma_2,\dots)$.
\hfill $\qed$

\medskip
The rotation strategy we used to count the marked Dyck paths is known for proving Narayana's formula. 

\begin{remark}
For $\gamma=(1,1,\dots)$ and $n\ge 1$, the sequence
\[ a_n = \sum_{k=1}^n \binom{n}{k-1}c_{n,k} = \sum_{k=1}^n k\,d_{n,k} \]
gives $1, 3, 10, 35, 126, 462, 1716, 6435, 24310, 92378,\dots$, cf.~\oeis{A001700}.
\end{remark}
% Dyck paths of semilength n having one of its peaks marked.

%%%%%%%%%%%%%%%%%%%%%%%%%%%%%%%%%%%%%%%%%%%%%%%%%%%%%
\section{Bijection between Dyck paths and polygon partitions}
\label{sec:Poly2Dyck}

In this section, we discuss a connection between colored Dyck paths and colored partitions of a convex polygon made by noncrossing diagonals. 

We start by pointing out that polygon partitions are easily converted to rooted trees. In fact, a convex $(n+2)$-gon partitioned by $k-1$ noncrossing diagonals can be bijectively mapped to a rooted tree with $n+1$ leaves, $k-1$ internal nodes, and having no node of outdegree 1. An example of this known bijection is illustrated in Figure~\ref{fig:polygon2tree}.

\begin{figure}[ht]
\begin{tikzpicture}
\begin{scope}
\node[regular polygon, regular polygon sides=6, minimum size=14ex] at (0,0) (A) {};
\draw[thick,gray] (A.corner 3) -- (A.corner 6);
\draw[thick,gray] (A.corner 4) -- (A.corner 6);
%outer edges
\node[regular polygon, regular polygon sides=6, minimum size=14.2ex, draw, thick] at (0,0) (B) {};
\draw[thick,gray,dash pattern=on 3pt off 1pt] (B.corner 4) -- (B.corner 5);
\end{scope}
\begin{scope}[xshift=110]
\node[regular polygon, regular polygon sides=6, minimum size=14ex] at (0,0) (A) {};
\node[regular polygon, regular polygon sides=6, minimum size=12.2ex, rotate=30] at (0,0) (M) {};
\draw[thick,gray!60] (A.corner 3) -- (A.corner 6);
\draw[thick,gray!60] (A.corner 4) -- (A.corner 6);
%tree
\draw[red1!50,thick] (M.corner 4) -- (M.corner 5);
\draw[red1!50,thick] (M.corner 4) -- (0.3,-0.5) -- (0,0) -- (M.corner 1); 
\draw[red1!50,thick] (0.3,-0.5) -- (M.corner 3); 
\draw[red1!50,thick] (0,0) -- (M.corner 6); 
\draw[red1!50,thick] (0,0) -- (M.corner 2); 
%outer edges
\node[regular polygon, regular polygon sides=6, minimum size=14.2ex, draw, thick] at (0,0) (B) {};
\draw[thick,gray,dash pattern=on 3pt off 1pt] (B.corner 4) -- (B.corner 5);
\foreach \m in {1,...,6} {\draw[fill] (M.corner \m) circle(0.03);}
\foreach \x/\y in {0/0,0.3/-0.5} {\draw[fill] (\x,\y) circle(0.03);}
\draw[fill,brown!50!black] (M.corner 4) circle(0.04);
\end{scope}
\begin{scope}[xshift=210]
\node[regular polygon, regular polygon sides=6, minimum size=12.2ex, rotate=30] at (0,0) (M) {};
\draw[red1] (M.corner 4) -- (M.corner 5);
\draw[red1] (M.corner 4) -- (0.3,-0.5) -- (0,0) -- (M.corner 1); 
\draw[red1] (0.3,-0.5) -- (M.corner 3); 
\draw[red1] (0,0) -- (M.corner 6); 
\draw[red1] (0,0) -- (M.corner 2); 
\foreach \m in {1,...,6} {\draw[fill] (M.corner \m) circle(0.03);}
\foreach \x/\y in {0/0,0.3/-0.5} {\draw[fill] (\x,\y) circle(0.03);}
\draw[fill,brown!50!black] (M.corner 4) circle(0.04);
\end{scope}
\begin{scope}[xshift=300]
\node at (0,0){
\begin{forest}
for tree={circle, fill, draw, s sep=17pt, inner sep=0.6pt, l=0}
[
    [ ]
    [ [ [] [] [] ] [] ]
]
\end{forest}
};
\end{scope}
\draw[->,cyan] (7.7,-1.03) to[out=0, in=250] (9.1,0) to[out=70, in=180] (10.2,0.8);
\end{tikzpicture}
\caption{Bijection between polygon partitions and rooted trees}
\label{fig:polygon2tree}
\end{figure}

\begin{thm}
Let $1\le k\le n$. For any coloring sequence $\gamma=(\gamma_j)_{j\in\mathbb{N}}$,
\begin{equation*}
 \binom{n+k}{k}d_{n,k}(\gamma) = \binom{n+1}{k}p_{n,k}(\gamma).
\end{equation*}
\end{thm}

\medskip
As before, we start by assuming $\gamma=(1,1,1,\dots)$. 

Every $n$-Dyck path with $k$ peaks has $n$ $\U$-steps, and it must have $k$ $\D$-steps that are each part of a peak. We let $\DP_{n,k}^{(k\;\textup{steps})}$ denote the set of such paths where $k$ of these $n+k$ steps are marked. On the other hand, a rooted convex $(n+2)$-gon must have one side designated as the base. We let $\Pol_{n,k}^{(k)}$ be the set of partitions of a rooted convex $(n+2)$-gon made by $k-1$ noncrossing diagonals, where $k$ of the $n+1$ non-base sides of the polygon are marked.

The above bijection between rooted polygon partitions and rooted trees (depicted in Figure~\ref{fig:polygon2tree}) gives a bijection between $\Pol_{n,k}^{(k)}$ and the set $\T_{n,k}^{(k)}$ of rooted trees having $n+1$ leaves and $k$ nodes of outdegree greater than 1 (for a total of $n+k+1$ nodes), where $k$ of the leaves are marked. We will give an explicit bijective map $\varphi_{dp}: \DP_{n,k}^{(k\;\textup{steps})}\to \T_{n,k}^{(k)}$ along with an illustrating example.

\begin{enumerate}[\quad$\circ$]
\item Given a marked Dyck path $P$, break the path at each valley so that there are $k$ shorter lattice paths with one peak each. 
\begin{equation}\label{eq:markedDP}
\begin{tikzpicture}[baseline=0]
\node at (0,0) {\dyckpath{0,0}{10}{1,1,1,0,1,0,1,1,0,0,0,0,1,0,1,0,1,1,0,0}{0.35}};
\foreach \x/\y in {-3.25/-0.65,-2.9/-0.3,-2.55/0.05,-1.85/0.05,1.64/-0.65} {
\draw[red1, line width=1pt] (\x,\y) -- (\x+0.22,\y+0.22);
}
\draw[red1, line width=1pt] (-1.59,0.25) -- (-1.37,0.03);
\foreach \x in {-1.97,-1.27,0.83,1.53,2.23}{\draw[gray,thick,dotted] (\x,0.8) -- (\x,-0.8);}
\end{tikzpicture}
\end{equation}
\item For each segment $\U^a\D^{b}$, create a {\em primitive} rooted tree of the form 
\tikzstyle{tree}=[circle, draw, fill, inner sep=0pt, minimum width=2pt]
\begin{tikzpicture}[baseline=-14]
\draw (0,0) node[tree]{} -- (-0.5,-0.6) node[tree]{}; 
\draw (0,0) node[tree]{} -- (-0.3,-0.6) node[tree]{}; 
\draw (0.15,-0.6) node{$\dots$}; 
\draw (0,0) -- (0.5,-0.6) node[tree]{};
\node[right=1] at (0.1,-0.15) {\scriptsize $b$};
\end{tikzpicture}
with $a+1$ leaves and the weight $b$ (omitted if $b=1$) assigned to the right-most edge. Mark the leaves of any branches that correspond to a marked step on the Dyck path. This gives an ordered tuple of marked trees $(T_1,\dots,T_k)$ corresponding to the given Dyck path.

\begin{center}
\begin{tikzpicture}
\node at (0,0) {
\begin{forest}
for tree={circle, fill, draw, s sep=9pt, inner sep=0.6pt, l=0}
[
	[\reddot, fill=red1] [\reddot, fill=red1] [\reddot, fill=red1] []
]
\end{forest}
};
\node at (1.4,0) {
\begin{forest}
for tree={circle, fill, draw, s sep=12pt, inner sep=0.6pt, l=0}
[
	[\reddot, fill=red1] [\reddot, fill=red1]
]
\end{forest}
};
\node at (2.6,0) {
\begin{forest}
for tree={circle, fill, draw, s sep=10pt, inner sep=0.6pt, l=0}
[
	[] [] []
]
\end{forest}
};
\node at (3.8,0) {
\begin{forest}
for tree={circle, fill, draw, s sep=12pt, inner sep=0.6pt, l=0}
[ 
	[] []
]
\end{forest}
};
\node at (4.75,0) {
\begin{forest}
for tree={circle, fill, draw, s sep=12pt, inner sep=0.6pt, l=0}
[
	[\reddot, fill=red1] []
]
\end{forest}
};
\node at (5.85,0) {
\begin{forest}
for tree={circle, fill, draw, s sep=10pt, inner sep=0.6pt, l=0}
[
	[] [] []
]
\end{forest}
};
\node[right=4] at (2.6,0.1) {\scriptsize 4};
\node[right=4] at (5.85,0.1){\scriptsize 2};
\foreach \x in {0.87,1.92,3.3,4.25,5.22}{\draw[gray,thick,dotted] (\x,0.4) -- (\x,-0.4);}
\end{tikzpicture}
\end{center}

\item For every $j\in\{1,\dots,k\}$, let $\ell_j$ be the number of unmarked leaves in $T_j$ and let $b_j$ be the weight of its right-most edge. Going left to right, we proceed with the following merging process:

If $\ell_1 \geq b_1$, place the root of $T_2$ into $T_1$ on its $b_1$-st unmarked leaf from the right. We denote the merged tree by $T_{1,2}$ and declare the $b_1-1$ right-most leaves of $T_1$ as inactive. The tree $T_{1,2}$ now has $\ell_1-b_1+\ell_2$ unmarked active leaves, and we assign to it the weight $b_2$.
If $\ell_1 < b_1$, do nothing and move to $T_2$. 

We repeat this merging procedure until we reach $T_k$ and let $(F_{v_1},\dots,F_{v_{m}})$ be the resulting ordered forest of merged primitive trees.

\begin{center}
\begin{tikzpicture}
\begin{scope}
\node at (0,0) {
\begin{forest}
for tree={circle, fill, draw, s sep=9pt, inner sep=0.6pt, l=0}
[
	[\reddot, fill=red1] [\reddot, fill=red1] [\reddot, fill=red1] []
]
\end{forest}
};
\draw[->,cyan] (1.4,-0.4) to[out=-135, in=-45] (0.7,-0.4);
\node at (1.4,0) {
\begin{forest}
for tree={circle, fill, draw, s sep=12pt, inner sep=0.6pt, l=0}
[
	[\reddot, fill=red1] [\reddot, fill=red1]
]
\end{forest}
};
\node at (2.65,0) {
\begin{forest}
for tree={circle, fill, draw, s sep=10pt, inner sep=0.6pt, l=0}
[
	[] [] []
]
\end{forest}
};
\node at (3.8,0) {
\begin{forest}
for tree={circle, fill, draw, s sep=12pt, inner sep=0.6pt, l=0}
[ 
	[] []
]
\end{forest}
};
\draw[->,cyan] (4.8,-0.4) to[out=-135, in=-45] (4.12,-0.4);
\node at (4.8,0) {
\begin{forest}
for tree={circle, fill, draw, s sep=12pt, inner sep=0.6pt, l=0}
[
	[\reddot, fill=red1] []
]
\end{forest}
};
\draw[->,cyan] (5.85,-0.4) to[out=-135, in=-45] (5.15,-0.4);
\node at (5.85,0) {
\begin{forest}
for tree={circle, fill, draw, s sep=10pt, inner sep=0.6pt, l=0}
[
	[] [] []
]
\end{forest}
};
\node[right=4] at (2.65,0.1) {\scriptsize 4};
\node[right=4] at (5.85,0.1){\scriptsize 2};
\end{scope}

\begin{scope}[yshift=-50]
\node at (0.5,0) {
\begin{forest}
for tree={circle, fill, draw, s sep=9pt, inner sep=0.6pt, l=0}
[
	[\reddot, fill=red1] [\reddot, fill=red1] [\reddot, fill=red1] 
	[ [\reddot, fill=red1] [\reddot, fill=red1] ]
]
\end{forest}
};
\node at (2.7,0.25) {
\begin{forest}
for tree={circle, fill, draw, s sep=10pt, inner sep=0.6pt, l=0}
[
	[] [] []
]
\end{forest}
};
\node at (5.1,0) {
\begin{forest}
for tree={circle, fill, draw, s sep=10pt, inner sep=0.6pt, l=0}
[
	[] [ [\reddot, fill=red1] 
	[ [] [] [] ] 
	] 
]
\end{forest}
};
\node[below=28pt] at (0.5,0) {\small $F_{v_1} = T_{1,2}$};
\node[below=28pt] at (2.75,0) {\small $F_{v_2} = T_{3}$};
\node[below=28pt] at (5.1,0) {\small $F_{v_3} = T_{4,5,6}$};
\node[right=7] at (2.6,0.35){\scriptsize 4};
\node[right=10] at (5,-0.45){\scriptsize 2};
\end{scope}
\end{tikzpicture}
\end{center}

If $m=1$, this step gives a tree in $\T_{n,k}^{(k)}$ which we denote by $\varphi_{dp}(P)$.

\item If $m>1$, we define $\varphi_{dp}(P)$ by means of an additional merging procedure. Note that a completely merged tree of the form $T_{i,\dots,j}$ has $\ell_i-b_i+\dots+\ell_{j-1}-b_{j-1}+\ell_j$ unmarked active leaves and weight $b_j$. Moreover, for $j<k$, we must have 
\[ \ell_i-b_i+\dots+\ell_{j-1}-b_{j-1}+\ell_j < b_j, \]
hence $\ell_i+\dots+\ell_{j}+\delta_{i,\dots,j} = b_i+\dots+ b_j$ for some $\delta_{i,\dots,j}\ge 1$. Thus, if $F_{v_{m}}=T_{i_m,\dots,k}$, the information from the trees $F_{v_1},\dots,F_{v_{m-1}}$ implies
\[ \ell_1+\dots+\ell_{i_m-1} +(\delta_{v_1}+\dots+\delta_{v_{m-1}}) = b_1+\dots+b_{i_m-1}, \]
where $\delta_{v_i}$ is the difference between the weight of the tree $F_{v_i}$ and the number of its unmarked active leaves. By definition, $\delta_{v_i}\ge 1$ for every $i\in\{1,\dots,m-1\}$. In addition, since $\sum\limits_{i=1}^k \ell_i=n=\sum\limits_{i=1}^k b_i$, we have $\ell_{i_m}+\dots+\ell_k-(\delta_{v_1}+\dots+\delta_{v_{m-1}}) =  b_{i_m}+\dots+b_k$, and therefore, 
\[ \ell_{i_m}-b_{i_m}+\dots+\ell_{k-1}-b_{k-1}+\ell_k = \delta_{v_1}+\dots+\delta_{v_{m-1}} + b_k. \]
This means that $F_{v_m}$ has exactly $\delta_{v_1}+\dots+\delta_{v_{m-1}}+b_k$ unmarked active leaves. Finally, we construct the tree $\varphi_{dp}(P)$ by attaching the roots of the trees $F_{v_1},\dots,F_{v_{m-1}}$ to the unmarked active leaves of $F_{v_{m}}$, from right to left, according to the pattern:

\begin{center}
\begin{tikzpicture}
\foreach \x in {2,4,5}{\draw[fill] (\x,0) circle(1pt);}
\foreach \x in {2,4,5}{\draw[gray] (\x,0) circle(1.8pt);}
\foreach \x in {2,4,5}{\node[below=0pt, gray] at (\x,0) {$\vert$};}
\foreach \x/\z in {2.1/$F_{v_{m-1}}$,4.1/$F_{v_2}$,5.1/$F_{v_1}$}{\node[below=14pt] at (\x,0) {\scriptsize \z};}
\foreach \x in {1.5,2.5,3.3,4.5,5.5}{\node at (\x,0) {$\cdots$};}
\draw[<-] (6.4,0) -- (7,0) node[right=5pt]{\parbox{5.2em}{\small\center all unmarked active leaves of $F_{v_{m}}$}};
\draw [decorate,decoration={brace,amplitude=3pt,raise=3pt}]
 (5.2,0) -- (5.72,0) node [midway,yshift=14pt,black] {\scriptsize $b_k$};
\draw [decorate,decoration={brace,amplitude=3pt,raise=3pt}]
 (3.9,0) -- (4.7,0) node [midway,yshift=14pt,black] {\scriptsize $\delta_{v_1}$};
\draw [decorate,decoration={brace,amplitude=3pt,raise=3pt}]
 (1.9,0) -- (2.7,0) node [midway,yshift=14pt,black] {\scriptsize $\delta_{v_{m-2}}$};
\end{tikzpicture}
\end{center}

In other words, the root of $F_{v_1}$ is attached to the right-most unmarked leaf of $F_{v_{m}}$ that allows the last $b_k$ active leaves of $F_{v_{m}}$ to stay unchanged. The remaining trees are attached in a way that $F_{v_{i+1}}$ is at the $\delta_{v_i}$-th unmarked leaf of $F_{v_{m}}$ to the left of $F_{v_i}$. Note that there will be $\delta_{v_{m-1}}\!-1\ge 0$ unmarked leaves to the left of $F_{v_{m-1}}$.

\begin{center}
\begin{tikzpicture}
\node at (0,0) {
\begin{forest}
for tree={circle, fill, draw, s sep=9pt, inner sep=0.6pt, l=0}
[
	[\reddot, fill=red1] [\reddot, fill=red1] [\reddot, fill=red1] 
	[ [\reddot, fill=red1] [\reddot, fill=red1] ]
]
\end{forest}
};
\node at (2.4,0.25) {
\begin{forest}
for tree={circle, fill, draw, s sep=10pt, inner sep=0.6pt, l=0}
[
	[] [] []
]
\end{forest}
};
\node at (5,0) {
\begin{forest}
for tree={circle, fill, draw, s sep=10pt, inner sep=0.6pt, l=0}
[
	[] [ [\reddot, fill=red1] 
	[ [] [] [] ] 
	] 
]
\end{forest}
};
\draw[gray] (4.463,0.28) circle(0.06);
\draw[gray] (4.685,-0.805) circle(0.06);
\node[right=7] at (2.3,0.35){\scriptsize 4};
\node[right=8] at (5,-0.45){\scriptsize 2};
\draw[->,cyan] (2.4,-0.3) to[out=-30, in=230] (4.4,0.15);
\draw[->,cyan] (0,-0.8) to[out=-30, in=210] (4.6,-0.9);
\end{tikzpicture}
\end{center}

For the marked Dyck path given in \eqref{eq:markedDP}, the above construction leads to the following marked tree and corresponding polygon partition:
\begin{center}
\begin{tikzpicture}
\begin{scope}
\node at (0,-0.2) {
\begin{forest}
for tree={circle, fill, draw, s sep=10pt, inner sep=0.6pt, l=0}
[
	[ [][][] ]
	[ [\reddot, fill=red1]
	  [ [[\reddot, fill=red1][\reddot, fill=red1][\reddot, fill=red1]
	  [ [\reddot, fill=red1][\reddot, fill=red1] ]] [] [] 
	] ]
]
\end{forest}
};
\end{scope}
\node at (2.4,0) {\large $\leftrightsquigarrow$};
\begin{scope}[xshift=140]
\node[regular polygon, regular polygon sides=12, minimum size=15ex] at (0,0) (A) {};
\draw[thick, gray] (A.corner 3) -- (A.corner 5);
\draw[thick, gray] (A.corner 5) -- (A.corner 12);
\draw[thick, gray] (A.corner 12) -- (A.corner 7);
\draw[thick, gray] (A.corner 7) -- (A.corner 11);
\draw[thick, gray] (A.corner 11) -- (A.corner 8);
%outer edges
\node[regular polygon, regular polygon sides=12, minimum size=15.2ex, draw, thick] at (0,0) (B) {};
\draw[thick,gray,dash pattern=on 3pt off 1pt] (B.corner 7) -- (B.corner 8);
\draw[very thick, red1] (B.corner 1) -- (B.corner 2);
\draw[very thick, red1] (B.corner 2) -- (B.corner 3);
\draw[very thick, red1] (B.corner 3) -- (B.corner 4);
\draw[very thick, red1] (B.corner 4) -- (B.corner 5);
\draw[very thick, red1] (B.corner 11) -- (B.corner 12);
\draw[very thick, red1] (B.corner 12) -- (B.corner 1);
\end{scope}
\end{tikzpicture}
\end{center}
\end{enumerate}

\medskip
Conversely, there is an algorithm to label and decompose each element of $\T_{n,k}^{(k)}$ into an ordered forest of primitive rooted trees that leads to an element of $\DP_{n,k}^{(k\;\textup{steps})}$. 

Suppose we are given $T\in\T_{n,k}^{(k)}$. Such a tree consists of $k$ primitive subtrees, $k-1$ internal nodes, $n+k$ edges, and $n+1-k$ unmarked leaves.
\begin{enumerate}[\quad$\circ$]
\item 
Starting at the root of $T$, and going around clockwise, denote its primitive components by $T_1,\dots,T_k$. For $i\in\{1,\dots,k-1\}$, label the right-most edge of $T_i$ with the number $\lambda_i\ge 1$ of unmarked nodes needed to arrive at the root of $T_{i+1}$. This gives a sequence of labels $\lambda_1,\dots,\lambda_{k-1}$ with $\lambda_1+\dots+\lambda_{k-1}\le n$. We let $\lambda_k=n-(\lambda_1+\dots+\lambda_{k-1})$.

\begin{center}
\begin{tikzpicture}
\node at (0,-0.2) {
\begin{forest}
for tree={circle, fill, draw, s sep=12pt, inner sep=0.8pt, l=0}
[
	[ [][][] ]
	[ [\reddot, fill=red1]
	  [ [[\reddot, fill=red1][\reddot, fill=red1][\reddot, fill=red1]
	  [ [\reddot, fill=red1][\reddot, fill=red1] ]] [] [] 
	] ]
]
\end{forest}
};
\foreach \x/\y in {0.28/1.05,0.73/0.5,0.8/-0.65,1.33/-1.25} {\node at (\x,\y) {\tiny 1};}
\foreach \x/\y in {0.95/-0.12,-0.65/0.5} {\node[right=0pt] at (\x,\y) {\tiny 3};}
\end{tikzpicture}
\end{center}
Let $a_i+1$ be the number of edges of $T_i$, so $a_1+\dots+a_k=n$. Let $d_i=a_i-\lambda_i$.
\item If $\sum\limits_{i=1}^j d_i\ge 0$ for every $j\in\{1,\dots,k-1\}$, then we let $\varphi_{dp}^{-1}(T)$ be the marked Dyck path $\U^{a_1}\D^{\lambda_1} \cdots \U^{a_k}\D^{\lambda_k}$, where a subpath $\U^{a_i}\D$ is marked at each step corresponding to the marked leaves of $T_i$.
Otherwise, if the above path does not satisfy the Dyck path condition, we let $\ell>1$ be the smallest index such that
\[ \sum_{i=\ell}^j d_i\ge 0 \;\text{ for every } j\in\{\ell,\dots,k-1\}. \]
In this case, $\varphi_{dp}^{-1}(T)$ is defined to be the Dyck path 
\begin{equation}\label{eq:DyckWord}
  \U^{a_\ell}\D^{\lambda_\ell} \cdots \U^{a_k}\D^{\lambda_k+1}\,\U^{a_1}\D^{\lambda_1} \cdots \U^{a_{\ell-1}}\D^{\lambda_{\ell-1}-1} 
\end{equation}
associated with the sequence $(T_\ell,\dots,T_k,T_1,\dots,T_{\ell-1})$, and marked accordingly. In our example, we have $n=10$, $k=6$, $\ell=4$, and the Dyck path is 
\begin{center}
\tikz[baseline=0]{\node at (0,0){$(\U\U\U\D)(\U\D)(\U\U\D\D\D\D)(\U\D)(\U\D)(\U\U\D\D)$.};
\foreach \x in {-3.48,-3.22,-2.95,-2.12,-1.85,1.51}{\draw[red1,very thick] (\x,-0.18) -- (\x+0.2,-0.18);}}
\end{center}
To verify that \eqref{eq:DyckWord} is indeed a Dyck path, note that by definition, $d_{\ell-1}<0$. If there is a $J<\ell-1$ such that $d_{J+1}+\dots+d_{\ell-1}<0$ but $d_J+d_{J+1}+\dots+d_{\ell-1}\ge 0$, then 
\[ d_J \ge -(d_{J+1}+\dots+d_{\ell-1}) > 0, \]
and so $\sum\limits_{i=J}^j d_i\ge 0$ for every $j\in\{J,\dots,k-1\}$. This contradicts the minimality of $\ell$.

Therefore, $\sum\limits_{i=j}^{\ell-1} d_i < 0$ for every $j\in\{1,\dots,\ell-1\}$. By definition, $\sum\limits_{i=1}^{k} d_i=0$, so if we let $v=\sum\limits_{i=\ell}^{k} d_i$, then $v>0$ and $v+\sum\limits_{i=1}^jd_i>0$ for every $j\in\{1,\dots,\ell-2\}$. In addition, note that since $v-1\ge 0$, we have $(a_\ell-\lambda_\ell)+\dots+(a_k-(\lambda_k+1))\ge 0$.
\end{enumerate}

Under our bijection, a primitive subtree with $a+1$ leaves ($(a+2)$-gon on the polygon) corresponds to a primitive block $\U^a\D$ on the Dyck path. Therefore, our proof applies verbatim to the case of an arbitrary coloring sequence $\gamma = (\gamma_1, \gamma_2,\dots)$. \hfill $\qed$

\begin{remark}
For $\gamma=(1,1,\dots)$ and $n\ge1$, the sequence
\[ a_n = \sum_{k=1}^n  \binom{n+k}{k}d_{n,k} = \sum_{k=1}^n \binom{n+1}{k}p_{n,k} \]
gives $2, 9, 54, 375, 2848, 22981, 193742, 1688427, 15101778, 137930199,\dots$ \oeis{A368178}.
\end{remark}
% Rooted trees with $n+1$ leaves, no node of outdegree 1, and having as many leaves marked as the number of nodes of outdegree greater than 1.
% For example, if n=2, there is one tree [ [][][] ] with only one node of outdegree>1 (the root). This tree leads to 3 marked trees. The tree [ [] [[][]] ]  has 2 nodes of outdegree>1, so it gives binomial(3,2) = 3 marked trees. Similarly, the tree [ [[][]] [] ] gives 3 more marked trees for a total of 9.

%%%%%%%%%%%%%%%%%%%%%%%%%%%%%%%%%%%%%%%%%%%%%%%%%%%%%
\section{Bijection between polygon dissections and compositions}
\label{sec:Poly2Comp}

We now aim at providing a combinatorial proof for the last identity in \eqref{eq:mainIdentities}:

\begin{thm}
Let $1\le k\le n$. For any coloring sequence $\gamma=(\gamma_j)_{j\in\mathbb{N}}$,
\begin{equation*}
 (n+1)p_{n,k}(\gamma) = \binom{n+k}{k}c_{n,k}(\gamma).
\end{equation*}
\end{thm}

As in the previous two sections, our bijection will be part preserving, so without loss of generality we assume $\gamma=(1,1,1,\dots)$. The general case follows verbatim. 

To prove the above identity, we will replace polygon dissections with rooted trees (as done in the previous section, see Figure~\ref{fig:polygon2tree}), and it is more convenient to think of a composition as a binary word starting with $\sf 1$ and having no consecutive $\sf 1$'s. More precisely, a composition of $n$ with $k$ parts, say $(j_1, \dotsc, j_k)$, can be represented by the binary word $w_{j_1}\cdots w_{j_k}$ of length $n+k$, where $w_{j_i}$ consists of a one followed by $j_i$ zeros. For example, the composition $(1,3,2,4)$ corresponds to the binary word $\sf 10100010010000$.

With this in mind, we let $\T_{n,k}^{(1)}$ be the set of rooted trees having $n+1$ leaves and $k$ nodes of outdegree greater than one, where one of the leaves is marked. Moreover, we let $\B_{n,k}^{(k)}$ be the set of binary words of length $n+k$ starting with $1$, having $n$ $\sf 0$'s and $k$ $\sf 1$'s, no two consecutive $\sf 1$'s, and where $k$ of the letters are marked. Clearly, 
\[ \abs{\T_{n,k}^{(1)}}=(n+1)p_{n,k} \;\text{ and }\; \abs{\B_{n,k}^{(k)}}=\binom{n+k}{k}c_{n,k}. \]
Our goal is to give a bijective map $\varphi_{bt}: \B_{n,k}^{(k)}\to \T_{n,k}^{(1)}$ along with illustrating examples.

\begin{enumerate}[\quad$\circ$]
\item Every binary word in $\B_{n,k}^{(k)}$ is of the form $w=w_{j_1}\cdots w_{j_k}$, where $w_j = {\sf 10}^{\,j}$ (a one followed by $j$ zeros) and $k$ of the letters in $w$ are marked. Split $w$ into its $k$ primitive components $(w_{j_1},\dots,w_{j_k})$, and for every component $w_{j_i}$ make a primitive rooted tree $T_i$ with $j_i+1$ leaves. Mark the leaves that correspond to marked letters in the binary word. For example, the word $\sf 10\rcol{1}0100\rcol{10010}10$ in $\B_{8,6}^{(6)}$ leads to the decomposition
\begin{equation}\label{eq:bin2tree}
\begin{tikzpicture}[baseline=0pt]
\node[above=14pt] at (0,0) {{\sf 10}};
\node at (0,0) {
\begin{forest}
for tree={circle, fill, draw, s sep=12pt, inner sep=0.6pt, l=0}
[ [][] ]
\end{forest}
};
\node[above=14pt] at (1,0) {{\sf \rcol{1}0}};
\node at (1,0){
\begin{forest}
for tree={circle, fill, draw, s sep=12pt, inner sep=0.6pt, l=0}
[
    [\reddot, fill=red1]
    []
]
\end{forest}
};
\node[above=14pt] at (2.15,0) {{\sf 100}};
\node at (2.15,0){
\begin{forest}
for tree={circle, fill, draw, s sep=10pt, inner sep=0.6pt, l=0}
[
    [][][]
]
\end{forest}
};
\node[above=14pt] at (3.5,0) {{\sf \rcol{100}}};
\node at (3.5,0){
\begin{forest}
for tree={circle, fill, draw, s sep=10pt, inner sep=0.6pt, l=0}
[
    [\reddot, fill=red1][\reddot, fill=red1][\reddot, fill=red1] 
]
\end{forest}
};
\node[above=14pt] at (4.65,0) {{\sf \rcol{10}}};
\node at (4.65,0){
\begin{forest}
for tree={circle, fill, draw, s sep=10pt, inner sep=0.6pt, l=0}
[
    [\reddot, fill=red1][\reddot, fill=red1]
]
\end{forest}
};
\node[above=14pt] at (5.55,0) {{\sf 10}};
\node at (5.55,0){
\begin{forest}
for tree={circle, fill, draw, s sep=12pt, inner sep=0.6pt, l=0}
[
    [][]
]
\end{forest}
};
\end{tikzpicture}
\end{equation}
\item Combine the sequence of primitive trees $(T_1,\dots,T_k)$ by shifting each component to the left (starting with $T_2$) so that its root is placed on the right-most marked leaf of the tree to its immediate left. If there is no marked leaf to place $T_i$, move to $T_{i+1}$ and continue the process until you reach $T_k$. This step gives an ordered forest $(F_1,\dots,F_m)$ of merged primitive trees with a total of $m$ marked leaves, all of them on $F_m$.

\begin{center}
\begin{tikzpicture}
\begin{scope}
\node at (0,0) {
\begin{forest}
for tree={circle, fill, draw, s sep=12pt, inner sep=0.6pt, l=0}
[ [][] ]
\end{forest}
};
\node at (1,0){
\begin{forest}
for tree={circle, fill, draw, s sep=12pt, inner sep=0.6pt, l=0}
[
    [\reddot, fill=red1]
    []
]
\end{forest}
};
\node at (2.15,0){
\begin{forest}
for tree={circle, fill, draw, s sep=10pt, inner sep=0.6pt, l=0}
[
    [][][]
]
\end{forest}
};
\node at (3.5,0){
\begin{forest}
for tree={circle, fill, draw, s sep=10pt, inner sep=0.6pt, l=0}
[
    [\reddot, fill=red1][\reddot, fill=red1][\reddot, fill=red1] 
]
\end{forest}
};
\node at (4.65,0){
\begin{forest}
for tree={circle, fill, draw, s sep=10pt, inner sep=0.6pt, l=0}
[
    [\reddot, fill=red1][\reddot, fill=red1]
]
\end{forest}
};
\node at (5.55,0){
\begin{forest}
for tree={circle, fill, draw, s sep=12pt, inner sep=0.6pt, l=0}
[
    [][]
]
\end{forest}
};
\draw[->,cyan] (2.15,-0.4) to[out=225, in=-45] (0.8,-0.4);
\draw[->,cyan] (4.65,-0.4) to[out=230, in=-50] (4,-0.4);
\draw[->,cyan] (5.55,-0.4) to[out=230, in=-50] (4.95,-0.4);
\end{scope}
\begin{scope}[xshift=15,yshift=-40]
\node at (0,0) {
\begin{forest}
for tree={circle, fill, draw, s sep=12pt, inner sep=0.6pt, l=0}
[ [][] ]
\end{forest}
};
\node at (1.7,-0.28){
\begin{forest}
for tree={circle, fill, draw, s sep=12pt, inner sep=0.6pt, l=0}
[
    [\reddot, fill=red1
      [][][]
    ]
    []
]
\end{forest}
};
\node at (4,-0.56){
\begin{forest}
for tree={circle, fill, draw, s sep=10pt, inner sep=0.6pt, l=0}
[
    [\reddot, fill=red1][\reddot, fill=red1]
    [\reddot, fill=red1
      	[\reddot, fill=red1]
      	[\reddot, fill=red1
      		[][]
      	]
    ] 
]
\end{forest}
};
\end{scope}
\draw[->,thick,gray] (6.4,-0.3) to[out=-20, in=20] (6.4,-1.8);
\end{tikzpicture}
\end{center}

\item If $m=1$, let $\varphi_{bt}(w)$ be the tree obtained from $F_1$ by removing the marks from its internal nodes. If $m>1$, attach the trees $F_1,\dots,F_{m-1}$ (from right to left) to the $m-1$ left-most marked leaves on $F_m$ and remove the internal marks from the merged tree. The resulting tree will be $\varphi_{bt}(w)\in \T_{n,k}^{(1)}$.

\begin{center}
\begin{tikzpicture}
\begin{scope}
\node at (0,0) {
\begin{forest}
for tree={circle, fill, draw, s sep=12pt, inner sep=0.6pt, l=0}
[ [][] ]
\end{forest}
};
\node at (1.7,-0.28){
\begin{forest}
for tree={circle, fill, draw, s sep=12pt, inner sep=0.6pt, l=0}
[
    [\reddot, fill=red1
      [][][]
    ]
    []
]
\end{forest}
};
\node at (4,-0.56){
\begin{forest}
for tree={circle, fill, draw, s sep=10pt, inner sep=0.6pt, l=0}
[
    [\reddot, fill=red1][\reddot, fill=red1]
    [\reddot, fill=red1
      	[\reddot, fill=red1]
      	[\reddot, fill=red1
      		[][]
      	]
    ] 
]
\end{forest}
};
\draw[->,cyan] (2,-0.4) to[out=-50, in=230] (3.28,-0.4);
\draw[->,cyan] (0,-0.4) to[out=-70, in=250] (3.76,-0.4);
\end{scope}
\node at (5.7,-0.5) {\large $\leadsto$};
\begin{scope}[xshift=230]
\node at (0,-0.6){
\begin{forest}
for tree={circle, fill, draw, s sep=10pt, inner sep=0.6pt, l=0}
[
    [ [ [][][] ][] ][ [][] ]
    [
      	[\reddot, fill=red1]
      	[ [][] ]
    ] 
]
\end{forest}
};
\end{scope}
\end{tikzpicture}
\end{center}
\end{enumerate}

The inverse map is straightforward. Let $T\in\T_{n,k}^{(1)}$. Such a tree consists of $k$ primitive subtrees, $k-1$ internal nodes, and one marked leaf.
\begin{enumerate}[\quad$\circ$]
\item Starting at the root of $T$, and going around clockwise, denote its primitive components by $T_1,\dots,T_k$. Let $T_m$ be the subtree with the marked leaf. 
\item If the nodes of $T_m$ to the right of the marked leaf have a total of $j$ descendants, consider the ordered forest $(T_{m+j+1},\dots,T_k,T_1,\dots,T_{m+j})$ and mark the leaves that correspond to internal nodes of $T$.

\begin{center}
\begin{tikzpicture}
\begin{scope}[yshift=55]
\node at (0,0){
\begin{forest}
for tree={circle, fill, draw, s sep=12pt, inner sep=0.6pt, l=0}
[
    [ [ [][][] ][] ][ [][] ]
    [
      	[\reddot, fill=red1]
      	[ [][] ]
    ] 
]
\end{forest}
};
\node at (0.8,0.7) {\tiny $T_1$};
\node at (1.5,0.1) {\tiny $T_2$};
\node at (1.75,-0.5) {\tiny $T_3$};
\node at (0.48,0.1) {\tiny $T_4$};
\node at (-1.25,0.1) {\tiny $T_5$};
\node at (-1.7,-0.5) {\tiny $T_6$};
\node at (2.7,0) {\large $\leadsto$};
\node at (5,0) {\small $(T_4,T_5,T_6,T_1,T_2,T_3)$};
\end{scope}
\begin{scope}
\node at (0,0) {
\begin{forest}
for tree={circle, fill, draw, s sep=12pt, inner sep=0.6pt, l=0}
[ [][] ]
\end{forest}
};
\node at (1,0){
\begin{forest}
for tree={circle, fill, draw, s sep=12pt, inner sep=0.6pt, l=0}
[
    [\reddot, fill=red1]
    []
]
\end{forest}
};
\node at (2.15,0){
\begin{forest}
for tree={circle, fill, draw, s sep=10pt, inner sep=0.6pt, l=0}
[
    [][][]
]
\end{forest}
};
\node at (3.5,0){
\begin{forest}
for tree={circle, fill, draw, s sep=10pt, inner sep=0.6pt, l=0}
[
    [\reddot, fill=red1][\reddot, fill=red1][\reddot, fill=red1] 
]
\end{forest}
};
\node at (4.65,0){
\begin{forest}
for tree={circle, fill, draw, s sep=10pt, inner sep=0.6pt, l=0}
[
    [\reddot, fill=red1][\reddot, fill=red1]
]
\end{forest}
};
\node at (5.55,0){
\begin{forest}
for tree={circle, fill, draw, s sep=12pt, inner sep=0.6pt, l=0}
[
    [][]
]
\end{forest}
};
\end{scope}
\draw[->,thick,gray] (7,1.95) to[out=0, in=0] (7,0);
\end{tikzpicture}
\end{center}
\item For a primitive tree with $j+1$ leaves, assign the word ${\sf 10}^{\,j}$ and mark the letters that correspond to the marked leaves of the tree (like in \eqref{eq:bin2tree}). Finally, $\varphi_{bt}^{-1}(T)$ is obtained by concatenating the $k$ words associated with $(T_{m+j+1},\dots,T_k,T_1,\dots,T_{m+j})$ into a single binary word. \hfill\qed
\end{enumerate}

\begin{remark}
For $\gamma=(1,1,\dots)$ and $n\ge1$, the sequence
\[ a_n = \sum_{k=1}^n (n+1)p_{n,k} = \sum_{k=1}^n \binom{n+k}{k}c_{n,k} \]
gives $2, 9, 44, 225, 1182, 6321, 34232, 187137, 1030490, 5707449,\dots$, cf.~\oeis{A176479}.
\end{remark}
% Rooted ordered trees with n+1 leaves, no node of outdegree 1, and having one of its leaves marked.

%%%%%%%%%%%%%%%%%%%%%%%%%%%%%%%%%%%%%%%%%%%%%%%%%%%%%
%%%%%%%%%%%%%%%%%%%%%%%%%%%%%%%%%%%%%%%%%%%%%%%%%%%%%

\end{document}